\documentstyle{article}

\def\b{\beta}

\def\G{\Gamma}

\def\t{\tau}

\def\th{\theta}
\def\l{\lambda}
\def\L{\Lambda}

\def\so{\underline}

\def\f{\rightarrow}
\def\q{\forall}

\def\v{\vdash}

\def\p{\succ}

\begin{document}

{\Large \bf Un R\'esultat de Compl\'etude pour les Types $\q^+$ \\du
Syst\`eme ${\cal F}$} \\

{\bf Samir FARKH et Karim NOUR}\\

LAMA - Equipe de Logique -
Universit\'e de Savoie -
73376 Le Bourget du Lac.\\
E-mail sfarkh,knour@univ-savoie.fr\\

{\bf R\'esum\'e.} Nous pr\'esentons dans cette note un r\'esultat de
compl\'etude pour les types
\`a quantificateurs positifs du syst\`eme ${\cal F}$ de J.-Y. Girard. Ce
r\'esultat g\'en\'eralise
un th\'eor\`eme de R. Labib-Sami (voir [3]). \\

{\Large \bf A Completeness Result for the $\q^+$ Types \\of System ${\cal
F}$} \\

{\bf Abstract.} We presente in this note a completeness result for the
types with positive
quantifiers of the J.-Y. Girard type system ${\cal F}$. This result
generalizes a theorem of
R. Labib-Sami (see [3]).\\

{\bf Abridged English Version.} The ${\cal F}$ type system have been
introduced by J.-Y. Girard
(see [1]). This system is based on the second order intuitionistic
propositional calculus, and
thus it gives the possibility to quantify on types. In addition to
strong normalisation theorem
which certifies the termination of programs, the system ${\cal F}$ has two
more properties : \\
It allows to write programs for all the fonctions whose termination can be
proved in the Peano's
second order arithmetics. \\
- It allows to define all the usual data types : booleans, integers,
lists,
etc. \\
 The semantics of the system ${\cal F}$ proposed by J.-Y. Girard and J.-L.
Krivine (see [1] and
[2]) consists in associating to each type $A$ a set of $\l$-terms $|A|$,
in
order to obtain the
following result: if a $\l$-term $t$ is of type $A$, then it belongs to
the
set $|A|$. This
result is known as the adequation lemma, and it allows to prove the strong
normalisation of system
${\cal F}$ and the unicity of data representation.\\
The converse of the adequation lemma (a completeness result) is not true.
The difficulty comes from
the interpretation of second order quantifier. R. Labib-Sami proved a
completeness result for the
types with positive quantifiers and for a semantic based on the sets
stable
by the $\b
\eta$-equivalence (see [3]). In this note we prove a refined result by
indicating that
week-head-expansion suffices. Hence, it presents a generalisation of R.
Labib-Sami's result. In
the end of this note, we deduce some consequences among which the fact
that
for each type with
positive quantifiers $A$, the set $|A|$ is stable by $\b$-equivalence.

 \section {Notations et d\'efinitions}

On d\'esignera par $\L$ l'ensemble des termes du $\l$-calcul pur, dits
aussi $\l$-termes.
\'Etant donn\'es des $\l$-termes $t, u, u_1,..., u_n$, l'application de
$t$
\`a $u$ sera
not\'ee $(t)u$, et $(... ((t)u_1) ...)u_n$ sera not\'e $(t)u_1...u_n$. Si
$t$ est un $\l$-terme, on
d\'esigne par $Fv (t)$ l'ensemble de ses variables libres. On note par
$\f_{\b}$ la $\b$-r\'eduction,
et par $\simeq\sb{\b}$ la $\b$-\'equivalence. Un $\l$-terme $t$ soit 
poss\`ede un redex de t\^ete
faible [i.e. $t = (\l x u) v  v_1 ... v_m$, le redex de t\^ete faible est
$(\l x u ) v$], soit 
est en forme normale de t\^ete faible [i.e. $t= (x ) v_1 ... v_m$ ou $t =
\l x v$]. La notation $u
\p_f v$ signifie que $v$ est obtenu \`a partir de $u$ par r\'eduction de
t\^ete faible.\\

Nous utilisons comme syst\`eme de typage le syst\`eme ${\cal F}$ de J.-Y.
Girard. Les types de ce
syst\`eme sont les formules construites \`a l'aide d'un ensemble
d\'enombrable de variables
propositionnelles $X, Y$,..., et deux connecteurs $\f$ et $\q$. \'Etant
donn\'es
 un $\l$-terme
$t$, un type $A$, et un contexte $\G = \{x_1 : A_1,..., x_n : A_n \}$, on
d\'efinit au moyen des
r\`egles suivantes la notion ``$t$ est typable de type $A$ dans le
contexte
$\G$''. Cette
notion est not\'ee $\G \v_{\cal F} t : A$.\\
$(1)$ $\G \v_{\cal F} x_i : A_i$  $ (1 \leq i \leq n)$.\\
$(2)$ Si $\G, x : B \v_{\cal F} t : C$, alors $\G \v_{\cal F} \l xt : B \f
C$.\\
$(3)$ Si $\G \v_{\cal F} u : B\f C$, et $\G \v_{\cal F} v : B$, alors $\G
\v_{\cal F} (u)v : C$.\\
$(4)$ Si $\G \v_{\cal F} t : A$, et $X$ ne figure pas dans $\G$, alors $\G
\v_{\cal F} t : \q
XA$.\\
$(5)$ Si $\G \v_{\cal F} t : \q XA$, alors, pour tout type $C$, $\G
\v_{\cal F} t : A[C/X]$.\\

Il est facile de voir que : Si $\G\v_{\cal F} t : A$ et $\G \subseteq
\G'$, alors $\G'\v_{\cal F} t : A$. Et si $\G\v_{\cal F} t : A$, alors
$\G'\v_{\cal F} t : A$, o\`u $\G'$ est la restriction de $\G$ aux
d\'eclarations contenant les variables libres de $t$.\\

Le syst\`eme ${\cal F}$ poss\`ede les propri\'et\'es suivantes (voir [2]): \\

{\bf Th\'eor\`eme 1} 
{\it (i) Si $\G \v_{\cal F} t : A$, et $t\f_\beta t'$, alors $\G \v_{\cal
F} t' : A$.\\
(ii) Si $\G \v_{\cal F} t : A$, alors $t$ est fortement normalisable.}\\

Une partie $G$ de $\L$ est dite satur\'ee si, quels que soient les termes
$t$ et $u$, on a :
($u \in G$ et $t \p_f u) \Rightarrow t \in G$. Il est clair que
l'intersection d'un ensemble de
parties satur\'ees de $\L$ est satur\'ee. \'Etant donn\'ees deux parties
$G$ et $G'$ de $\L$,
on d\'efinit une partie de $\L$, not\'ee $G\f G'$, en posant : $ u\in (G\f
G') \Leftrightarrow
(u)t \in G'$ quel que soit $t\in G$. Si $G'$ est satur\'ee, alors $G\f G'$
est satur\'ee
pour toute partie $G \subset \L$. Une interpr\'etation $I$ est, par
d\'efinition, une application
$X \f  |X|_I$ de l'ensemble des variables de type dans l'ensemble des
parties satur\'ees de
$\L$. $X$ \'etant une variable de type, et $G$ une partie satur\'ee de
$\L$, on d\'efinit une
interpr\'etation $J = I [X \leftarrow G]$ en posant $|X|_J = G$, et $|Y|_J
= |Y|_I$ pour toute
variable $Y \neq X$. Pour chaque type $A$, sa valeur $|A|_I$ dans
l'interpr\'etation $I$ est une
partie satur\'ee  d\'efinie comme suit, par induction sur $A$ :\\
- Si $A$ est une variable de type, $|A|_I$ est d\'ej\`a d\'efinie ;\\
- $|A \f B|_I = |A|_I \f |B|_I$ ;\\
- $|\q XA|_I = \cap \{ |A|_{I [X \leftarrow G]}$ pour toute partie
satur\'ee $G \}$.\\

Il est facile de v\'erifier que : si $A, F$ sont deux types, $X$ une
variable, et $I$ une
interpr\'etation, alors $|A [F / X]|_I = |A|_{I [X \leftarrow |F|_I]}$.\\

Pour tout type $A$, on note $|A| = \cap \{|A|_I$ ; $I$
interpr\'etation$\}$. \\

{\bf Th\'eor\`eme 2} [lemme d'ad\'equation] {\it Soient $A$ un type, et
$t$
un $\l$-terme clos. \\
Si $\v_{\cal F} t : A$, alors $t \in |A|$.}

\section {Le r\'esultat de compl\'etude}

On d\'efinit de la fa\c con suivante les types \`a quantificateurs
positifs
(resp. \`a
quantificateurs n\'egatifs), not\'es en abr\'eg\'e  $\q^+$ (resp. $\q^-$)
:\\
- Une variable propositionnelle $X$ est $\q^+$ et $\q^-$ ;\\
- Si $A$ est $\q^+$ (resp. $\q^-$) et $B$ est $\q^-$ (resp. $\q^+$), alors
$B \rightarrow A$ est
$\q^+$ (resp. $\q^-$) ;\\
- Si $A$ est $\q^+$ et $X$ est libre dans $A$, alors $\q XA$ est $\q^+$.\\

On se propose de d\'emontrer le th\'eor\`eme suivant :\\

{\bf Th\'eor\`eme 3} {\it Soient A un type $\q^+$ du syst\`eme ${\cal F}$,
et $t$ un
$\l$-terme, alors : \\ $t \in |A|$ $\Leftrightarrow$ ($t \f_{\b} t'$ et
$\v_{\cal F} t' : A$).}\\

Pour la preuve, nous avons besoin de deux lemmes (lemme 1 et lemme 2).\\

{\bf Lemme 1} {\it Soient $I$ une interpr\'etation, et $t'$ un $\l$-terme
normal. Si
 $\G = x_1 : B_1,...,x_n : B_n \v_{\cal F} t' : A$, $t \simeq\sb{\b} t'$,
et  $u_i \in
|B_i|_I$ ($1 \leq i \leq n$), alors $t [u_1/x_1,...,u_n/x_n] \in
|A|_I$}.\\

 {\bf Preuve} : Par induction sur le typage. On consid\`ere la derni\`ere
r\`egle utilis\'ee.
\begin {itemize}
\item Si c'est la r\`egle (1), alors $t' = x_i$ ($1\leq i \leq n$) et $B_i
= A$. Comme
 $t \simeq\sb{\b} x_i$, alors $t \p_f x_i$, et $t [u_1/x_1,...,u_n/x_n]
\p_f u_i$. Or $u_i \in
|B_i|_I$, donc $t [u_1/x_1,...,u_n/x_n] \in |B_i|_I$, car $|B_i|_I$ est
une
partie satur\'ee.
\item Si c'est la r\`egle (2), alors $t' = \l xu'$, $A = B \f C$ et
 $\G, x : B \v_{\cal F} u' : A$. Comme $t \simeq\sb{\b} \l xu'$, alors $t
\p_f \l xu$ avec
 $u \simeq\sb{\b} u'$, et $t [u_1/x_1,...,u_n/x_n] \p_f \l xu
[u_1/x_1,...,u_n/x_n]$. Donc,
d'apr\`es l'hypoth\`ese d'induction, $u [u_1/x_1,...,u_n/x_n,v/x] \in
|C|_I$ pour tout $v \in |B|_I$.
D'autre part $(\l xu[u_1/x_1,...,u_n/x_n])v \p_f u
[u_1/x_1,...,u_n/x_n,v/x]$, donc $\l xu
[u_1/x_1,...,u_n/x_n] \in |B \f C|_I$, par cons\'equent $t
[u_1/x_1,...,u_n/x_n] \in |A|_I$.
\item Si c'est la r\`egle (3), comme $t'$ est normal, alors $t' = (u)v$,
 $\G \v_{\cal F} u : B\f A$ et  $\G \v_{\cal F} v : B$, avec
$u = (x_r)v'_1...v'_{m-1}$ et $v = v'_m$. Or $t \simeq\sb{\b}
(x_r)v'_1...v'_m$, donc
 $t \p_f (x_r)v_1...v_m$, avec $v_i \simeq\sb{\b} v'_i$ ($1\leq i \leq
m$),
d'o\`u, d'apr\`es
l'hypoth\`ese d'induction, \\
$(u_r)v_1 [u_1/x_1,...,u_n/x_n]...v_{m-1} [u_1/x_1,...,u_n/x_n] \in |B\f
A|_I$, et $v_m
[u_1/x_1,...,u_n/x_n] \in |B|_I$. Par cons\'equent
$(u_r)v_1 [u_1/x_1,...,u_n/x_n]...v_{m-1} [u_1/x_1,...,u_n/x_n]v_m
[u_1/x_1,...,u_n/x_n] \in |A|_I$,
et donc $t [u_1/x_1,...,u_n/x_n] \in |A|_I$.
\item Si c'est la r\`egle (4), alors  $\G \v_{\cal F} t' : B$ et $A = \q
XB$ avec $X$ ne figure pas
dans $\G$. Soit $G$ une partie satur\'ee et $J = I [X \leftarrow G]$. Par
hypoth\`ese d'induction,
$t [u_1/x_1,...,u_n/x_n] \in |B|_J$, et donc $t [u_1/x_1,...,u_n/x_n] \in
|A|_I$.
\item Si c'est la r\`egle (5), alors $\G \v_{\cal F} t' : \q XB$ et $A =
B
[C/X]$. Par hypoth\`ese
d'induction, $t [u_1/x_1,...,u_n/x_n] \in |\q XB|_I$, d'o\`u $t
[u_1/x_1,...,u_n/x_n] \in |B|_{I [X
\leftarrow |C|_I]} = |A|_I$. \hfill $\spadesuit$
\end {itemize}

Soient $\Omega = \{ x_i/i \in {\bf N} \}$ une \'enum\'eration d'un
ensemble
infini de variables du
$\l$-calcul, et $\{A_i/i \in {\bf N} \}$ une \'enum\'eration des types
$\q^-$ du syst\`eme ${\cal
F}$, o\`u chaque type $\q^-$ se r\'ep\`ete une infinit\'e de fois.
 On d\'efinit alors l'ensemble $\mit \G^- = \{x_i : A_i /i \in {\bf N}
\}$. Soit $u$ un
$\l$-terme, tel que $Fv (u) \subseteq \Omega$, on d\'efinit le contexte
$\mit \G^-_u$ comme
\'etant la restriction de $\mit \G^-$ aux d\'eclarations contenant les
variables de $Fv
(u)$. La notation $\mit \G^- \v_{\cal F} u : B$ exprime que $\mit \G^-_u
\v_{\cal F} u : B$. On
pose $\mit \G^- \v_{\cal F}^+ u : A$ ssi il existe un $\l$-terme $u'$, tel
que $u \f_{\b} u'$ et
$\mit \G^- \v_{\cal F} u' : A$. On d\'efinit ensuite une interpr\'etation
$\cal I$ en posant
$|X|_{\cal I} = \{ \t \in \L : \mit \G^-\v_{\cal F}^+\t : X \}$ pour toute
variable de type $X$.
Les parties $|X|_{\cal I}$ sont \'evidemment satur\'ees.\\

{\bf Lemme 2} {\it (i) Si $S$ est un type $\q^+$, et $\t \in |S|_{\cal
I}$,
alors $\mit \G^-
\v_{\cal F}^+ \t : S$.\\ (ii) Si $S$ est un type $\q^-$, et $\mit \G^-
\v_{\cal F}^+ \t : S$, alors
$\t \in |S|_{\cal I}$.}\\

{\bf Preuve} : Par induction simultan\'ee sur les types $\q^+$ et
$\q^-$. \\

\so{Preuve de (i)}
\begin {itemize}
\item Si $S$ est une variable, alors le r\'esultat d\'ecoule
imm\'ediatement de la d\'efinition de
${\cal I}$.

\item Si $S=\q XB$, o\`u $B$ est $\q^+$, alors soit $\t\in |\q XB|_{\cal
I}$, et soit $Y$ une
variable propositionnelle qui ne figure pas dans $\mit \G^-_\t$ et $B$. Donc
$\t\in |B|_{{\cal I} [X
\leftarrow |Y|_{\cal I}]} = |B[Y / X]|_{\cal I}$, d'o\`u par hypoth\`ese
d'induction $\mit \G^-
\v_{\cal F}^+ \t : B[Y / X]$, donc $\t \f_{\b} \t'$ et $\mit \G^-_{\t'}
\v_{\cal F} \t' :
B[Y / X]$. Comme $Fv (\t') \subseteq Fv (\t)$, alors par le choix de $Y$,
on d\'eduit que $\mit
\G^-_{\t'} \v_{\cal F} \t' : \q YB[Y / X]=\q XB$, et donc $\mit \G^- \v_{\cal
F}^+ \t : S$.

\item Si $S=B\f C$, o\`u $B$ est $\q^-$ et $C$ est $\q^+$, alors soit
$\t\in |B\f C|_{\cal I}$ et
soit $y$ une variable du $\l$-calcul telle que $y : B$ appartient \`a
$\mit
\G^-$. On a $y : B
\v_{\cal F} y : B$, donc d'apr\`es (ii) $y\in |B|_{\cal I}$, par suite
$(\t)y\in |C|_{\cal I}$, et
donc, d'apr\`es l'hypoth\`ese d'induction, $\mit \G^- \v_{\cal F}^+ (\t)y
:
C$. D'o\`u $(\t)y
\f_{\b} \t'$, et $\mit \G^-_{\t'} \v_{\cal F} \t': C$. Il en r\'esulte que
$(\t)y$ est
normalisable, et donc $\t$ est normalisable. La forme normale de $\t$ est
$(x)\t_1...\t_n$ $(n \geq 0)$ ou $\l x \th$.\\

{\so{\bf Cas 1}}: Si $\t\f_{\b} (x)\t_1...\t_n$, avec $n \geq 0$, alors
$(\t)y\f_{\b}(x)\t_1...\t_ny$. Comme $\mit
\G^-_{\t'} \v_{\cal F} \t' : C$, on aura $\mit
\G^-_{\t'} \v_{\cal F} (x)\t_1...\t_ny : C$.
 Or $Fv ((x)\t_1...\t_ny) \subseteq Fv (\t')$, donc
  $\mit \G^-_{(x)\t_1...\t_ny} \v_{\cal F} (x)\t_1...\t_ny : C$. Donc $x$
est
   d\'eclar\'ee de type $E$ dans le contexte $\mit
\G^-_{(x)\t_1...\t_ny}$. Comme $E$ est $\q^-$, alors $E = C_1 \f (... \f
(C_n \f (F \f G)...))$, $\mit
\G^-_{(x)\t_1...\t_ny} \v_{\cal F} \t_i : C_i$ ($1 \leq i \leq n$) et
$\mit
\G^-_{(x)\t_1...\t_ny} \v_{\cal F} y : F$. Or $y$ est d\'eclar\'ee de type
$B$,
 qui est $\q^-$, dans le contexte $\mit
\G^-_{(x)\t_1...\t_ny}$, donc on ne peut pas appliquer la r\`egle de
typage $(5)$
 \`a $y : B \v_{\cal F} y : B$.
  D'autre part, si on applique la r\`egle de typage $(4)$ \`a
  $y : B \v_{\cal F} y : B$, $F$ sera le type $B$ devant lequel on a
quantifi\'e
  sur une variable qui n'est pas libre dans $B$, ce qui contredit le fait
que 
  $F$ est $\q^+$. D'o\`u $F = B$. De la m\^eme fa\c con, comme $G$ est
$\q^-$
  et $C$ est $\q^+$, alors de
   $\mit \G^-_{(x)\t_1...\t_ny} \v_{\cal F} (x)\t_1...\t_ny : C$ et
    $\mit \G^-_{(x)\t_1...\t_ny} \v_{\cal F} (x)\t_1...\t_ny : G$, on
d\'eduit 
    que $G = C$. Par cons\'equent
     $\mit \G^-\v_{\cal F} (x)\t_1...\t_n : B\f C$, et donc $\mit
\G^-\v_{\cal F}^+ \t : S$.

{\so{\bf Cas 2}}: Si $\t\f_{\b}\l x \th$, alors, comme l'ensemble $\G^-$ contient une infinit\'e de
d\'eclarations pour chaque type $\q^-$, soit $y$ une variable
 d\'eclar\'ee de type $B$ dans  $\G^-$, n'appartenant pas
  \`a $Fv (\th)$. Alors $(\t)y \f_{\b}(\l x \th)y \f_{\b} \th[y/x]$.
   On a $\mit \G^-_{\t'} \v_{\cal F} \t' : C$, donc $\mit
\G^-_{\t'} \v_{\cal F}
\th[y/x] : C$, et $\mit \G^-_{\th[y/x]} \v_{\cal F} \th[y/x] : C$,
car
$Fv (\th[y/x])
\subseteq Fv (\t')$. Donc $\mit \G^-_{\l y \th[y/x]}
\v_{\cal F} \l y
\th[y/x] : B \f
C$. Comme $y$ n'est pas libre dans $\th$, alors $\l y \th [y/x] = \l x \th$, par cons\'equent $\mit
\G^- \v_{\cal F} \l x\th : S$ et $\mit \G^- \v_{\cal F}^+ \t : S$.
\end {itemize}

\so{Preuve de (ii)}
\begin {itemize}
\item Si $S$ est une variable, alors le r\'esultat d\'ecoule
imm\'ediatement de la d\'efinition de
${\cal I}$.
\item Si $S = B \rightarrow C$, o\`u $B$ est $\q^+$ et $C$ est $\q^-$,
alors supposons que $\mit
\G^-\v_{\cal F}^+ \t : B \rightarrow C$, donc $\t \f_{\b} \t'$ et $\mit
\G^-_{\t'}\v_{\cal F} \t' :
B \rightarrow C$. Si $u \in |B|_{\cal I}$, alors d'apr\`es (i), $\mit
\G^-\v_{\cal F}^+  u : B$,
donc $u \f_{\b} u'$ et $\mit \G^-_{u'}\v_{\cal F} u' : B$. D'o\`u $\mit
\G^-_{(\t')u'}\v_{\cal F}
(\t') u' : C$, et comme $(\t)u \f_{\b} (\t')u'$, alors  $\mit \G^-\v_{\cal
F}^+ (\t)u : C$. D'o\`u,
d'apr\`es l'hypoth\`ese d'induction, $(\t)u \in |C|_{\cal I}$. Par
cons\'equent $\t \in |B \f
C|_{\cal I}$. \hfill $\spadesuit$
\end {itemize}

{\bf Preuve du th\'eor\`eme 3}\\
$\Leftarrow$) D'apr\`es le lemme 1.\\
$\Rightarrow$) Si $t \in |A|$, alors $t \in |A|_{\cal I}$ pour toute
interpr\'etation ${\cal I}$
associ\'ee \`a un ensemble $\mit \G^-$. Or on peut supposer que $\mit
\G^-$
ne contient pas de
d\'eclarations de variables libres de $t$, d'o\`u, d'apr\`es le (i) du
lemme 2 et le fait que $Fv
(t') \subseteq Fv (t)$, on a $t \f_{\b} t'$ et $\v_{\cal F}  t' : A$.
\hfill $\spadesuit$\\

{\bf Corollaire 1} {\it Soient $A$ un type $\q^+$ du syst\`eme ${\cal F}$,
et $t$ un $\l$-terme. \\
(i) Si $t \in |A|$, alors $t$ est normalisable et $\b$-\'equivalent \`a un
terme clos.\\
(ii) $|A|$ est stable par $\b$-\'equivalence (i.e si $t \in |A|$ et $t
\simeq\sb{\b} t'$, alors
$t' \in |A|$).}\\
 
{\bf Preuve} : La partie (i) est une cons\'equence directe du th\'eor\`eme
3 et la partie (ii) se
d\'eduit du th\'eor\`eme 3 et du lemme 1. \hfill $\spadesuit$\\
   
On consid\`ere le syst\`eme de typage ${{\cal F}0}$ qui n'est autre que le
syst\`eme ${\cal F}$,
o\`u on remplace la r\`egle de typage (5) par la r\`egle :
$\G\v_{{\cal F}0} t : \q XA \Rightarrow \G\v_{{\cal F}0} t : A[Y/X] $ o\`u
$Y$ est une variable.\\

Il \'etait d\'emontr\'e dans [4] que si $A$ est un type $\q^+$ du
syst\`eme
${\cal F}$, et $t$ est
un $\l$-terme normal clos, alors $\v_{\cal F} t : A$ $\Rightarrow$
$\v_{{\cal F}0} t : A$. On a
donc le th\'eor\`eme suivant : \\

{\bf Th\'eor\`eme 4} {\it Soient $A$ un type $\q^+$ du syst\`eme ${\cal
F}$, et $t$ un
$\l$-terme, alors :\\ $t \in |A|$ $\Leftrightarrow$ ($t \f_{\b} t'$ et
$\v_{{\cal F}0} t' : A$).}\\

{\bf Remarques}. (1) Si on interpr\`ete les types par des ensembles
$\b$-satur\'es
 [un ensemble $G$ est dit $\b$-satur\'e si, quels que soient les termes
$t$
  et $u$, on a : $u \in G$ et $t \f_{\b} u$ $\Rightarrow$ $t \in G$],
alors
   le lemme 2 reste vrai ainsi que le th\'eor\`eme 3. On peut donc
d\'eduire que les interpr\'etations (satur\'ee et $\b$-satur\'ee)
d'un type $\q^+$ sont \'egales.\\
(2) Les types de donn\'ees peuvent \^etre d\'efinis
dans
le syst\`eme ${\cal F}$
par des types $\q^+$ clos. Par exemple : le type bool\'een est la formule
:
$Bool = \q X \{X\f
(X\f X)\}$, le type des entiers est la formule $Ent = \q X \{(X\f X)\f
(X\f
X)\}$, et le type des
listes d'\'el\'ements de type $Ent$ est donn\'ee par la formule $LEnt = \q
X \{(Ent \f (X\f X))\f
(X\f X)\}$. On peut donc d\'eduire du th\'eor\`eme 4 que : pour tout
type de donn\'ees $D$ du syst\`eme ${\cal F}$, $|D| = \{t\in \L$, $t
\f_{\b} t'$ et $\v_{\cal F} t'
: D \}$.\\
(3) Pour obtenir son r\'esultat, R. Labib-Sami a autoris\'e la
quantification sur
une variable m\^eme si elle ne figure pas dans le type. Dans ce cas le
th\'eor\`eme 3 n'est plus
vrai. Par exemple si on prend le type $A = \q X \{ (X \f \q YX)\f (X\f
X)\}$, alors on peut
remarquer facilement que le terme $I = \l xx \in |A|$ et $\not \v_{{\cal
F}} I : A$ (mais $I
\simeq\sb{\b \eta} I' = \l x \l y (x)y$ et $\v_{{\cal F}} I' : A$). \\[2cm]

{\bf Remerciement}. Nous remercions le rapporteur et C. Raffalli pour
leurs
 remarques.\\

\end{document}